\magnification=1200

\def \v {\vskip .3cm}
\def \vv {\vskip .15cm}

\centerline {\bf Singularities and bifurcations of}

\centerline {\bf 3-dimensional Poisson structures}

\v 

\centerline {by} 

\vv 

\centerline {J.P. Dufour}

\centerline {\it Universit\'e Montpellier II, Case 051, Pl. E. Bataillon,} 

\centerline {\it 34095 Montpellier Cedex 5, France}

\centerline {\it e-mail: dufourj@math.univ-montp2.fr} 

\v 

\centerline {and}

\vv 

\centerline {M. Zhitomirskii} 

\centerline {\it Department of Mathematics, Technion, 32000 Haifa, Israel} 

\centerline {\it e-mail: mzhi@techunix.technion.ac.il}

\v

\vv 

\centerline {\rm ABSTRACT} 

{\sl We  give a normal form for 
families of 3-dimensional Poisson structures. This 
allows us to classify singularities with nonzero 1-jet and 
typical bifurcations. The Appendix contains corollaries on 
classification of families of integrable 1-forms on $R^3$.}

\v

\centerline {\bf 1. Introduction and main results. } 

\v

Poisson structures are central objects in classical
mechanics and its quantization, at least on the mathematical level.
 The first extensively studied Poisson structures 
were the regular ones ([Li]), later appeared an interest
in the study of singularities ([We1]). 
In the 2-dimension case, the singularities of Poisson 
structures were classified by 
 V. Arnold [Ar].  
In this case the classification  is similar to the classification 
 of functions. The latter is not true beginning from the 3-dimensional 
case because of the Jacobi identity which starts to play 
an important role in all classification results.  

This paper is  devoted to the 3-dimensional case. 
We classify local families $P_\epsilon $ of 
Poisson structures on $R^3$ such that $P_0$ has 
a singular point $0\in R^3$ (i.e., $P_0(0)=0$) and  
$j^1_0P_0$ does not vanish. This allows us to classify 
individual germs of Poisson structures and 
 bifurcations which hold in generic 
1-parameter families.    

All objects are assumed to be of  the class 
$C^\infty $, all families are assumed to depend smoothly on 
parameters. Let $M$ be a smooth manifold and $N$ the ring of
smooth functions on $M$. Recall that a Poisson structure  on 
$M$ is a composition law  $(f,g)\mapsto\{ f,g\}$ on $N$ which endows $N$
with a Lie algebra structure and satisfies the condition   
$\quad\{fg,h\}= f\{g,h\}+g\{f,h\}$
for every $f$, $g$ and $h$ in $N.$ Such a composition law is called a
Poisson bracket. It is locally determined by the brackets $\{x_i,x_j\}$
where $(x_1,\dots ,x_n)$ are local coordinates. 
This Poisson structure can  also 
be viewed as the  
2-vector 
$P=\sum_{i<j} \{x_i,x_j\} \partial /\partial x_i\wedge\partial 
/\partial x_j$, then the Poisson bracket of two functions $f$ and $g$ 
is given by the relation $\{ f,g\}=P(df,dg).$  
By local equivalence of Poisson structures we always mean the 
equivalence with respect to the natural action of the group of local  
diffeomorphisms. Two local families $P_\epsilon $ and 
$\tilde P_\epsilon $ are called equivalent if there exists a 
family of diffeomorphisms $\phi _\epsilon $ such that $\phi _0(0)=0$  
and $(\phi _\epsilon )_*P_\epsilon = \tilde P_\epsilon $ for all 
small $\epsilon $. In this definition it is not required that 
$\phi _\epsilon (0) = 0$ as $\epsilon \ne 0$.

If $M$ is a 3-dimensional orientable manifold then any volume form 
$\Omega $ induces a 1-1 correspondence between Poisson structures and 
integrable 1-forms. This correspondence is realized 
by the isomorphism 
$b: Z\mapsto i_Z\Omega$ between $p$-vectors and $(3-p)$-forms. Therefore 
all results of this paper can be reformulated as classification results 
for local integrable 1-forms on $R^3$ defined up to multiplication by 
a nonvanishing function, see Appendix.

Our main result, allowing to analyse singularities and bifurcations, is a 
versal unfolding of any Poisson structure having nonzero 1-jet at the 
origin.  

\vskip .3cm 

{\bf Theorem 1.1.} 
{\it Let $P_\epsilon $ be a local family of Poisson structures 
on $R^3$ such that 
$P(0)=0$, \ $j^1_0P\ne 0$, where $P=P_0$. Assume that $j^1_0P$ is not 
isomorphic to the Poisson structure  
$$\{x,y\}=0, \ \{y,z\} = y, \ \{z,x\} = -x.\eqno (1.1)$$ 
Then $P_\epsilon $ is equivalent to a local family of the form 
$$\{x,y\}=z, \ \{y,z\} = U_\epsilon (x,y,z), \ \{z,x\} = V_\epsilon 
(x,y,z)\eqno (1.2)$$ and there exists a family of formal changes 
of the coordinates $x,y,z$, centered at the point $x=y=z=0$ for all 
$\epsilon $, reducing (1.2) to the form   
$$\{x ,y \}=z , \ \{y ,z \} = 
{\partial \hat f_\epsilon (x ,y )\over \partial x } + 
z {\partial \hat g_\epsilon (x ,y )\over \partial x }, \ 
\{z ,x \} = {\partial \hat f_\epsilon (x ,y )\over \partial y } + 
z {\partial \hat g_\epsilon (x ,y )\over \partial y },\eqno (1.3)$$
where the formal series $\hat f_\epsilon $ and $\hat g_\epsilon $ 
satisfy the relations  
$$\hat f_\epsilon (0) = \hat g_\epsilon (0) = 0, \ 
d\hat f_\epsilon \wedge d\hat g_\epsilon \equiv 0.\eqno (1.4)$$ }

\vskip .3cm

This normal form is, to some extend, a generalization of the Bogdanov-Takens 
formal normal form  
for families of vector fields on a plane (see [AI]). Given a family of vector fields  
$v_\epsilon (x,z) = \alpha _\epsilon (x,z){\partial \over \partial x} + 
\beta _\epsilon (x,z){\partial \over \partial z}$ on the $(x,z)$ plane, we 
can associate to it a family of Poisson structures on $R^3$ of the form 
$${\partial \over \partial y}\wedge v_\epsilon (x,z). \eqno (1.5)$$ 
It follows from the proof of Theorem 1.1 
in section 5 
that, in general, such a family is reducible to the normal form 
(1.3), where the series $\hat f_\epsilon (x ,y )$ and   
$\hat g_\epsilon (x ,y )$ do not depend on $y$. Therefore  
(1.3) takes  
the form  ${\partial \over \partial y}\wedge (z{\partial \over \partial x} + 
(A_\epsilon (x) + zB_\epsilon (x)){\partial \over \partial z})$   
which corresponds to the Bogdanov-Takens normal form.  

\vskip .3cm

Using the normal form (1.3)-(1.4) we distinguish the following singularity classes. 
In what follows $\hat f = \hat f_0$ and $\hat g = \hat g_0$, where 
$\hat f_\epsilon $ and 
$\hat g_\epsilon $ are the functional parameters in the normal form (1.3).

\vskip .2cm

1) A Poisson structure $P$ has a $V$ singularity at the origin if either 
$j^1_0P$ is isomorphic to the Poisson structure (1.1) or 
$j^1_0\hat g\ne 0$.

\vskip .2cm

2) A Poisson structure  has a $so(3)$ singularity at the origin if  
$j^2_0\hat f$ is R-equivalent to $x^2 + y^2$ (then $j^1_0\hat g = 0$ by (1.4)).

\vskip .2cm

3) A Poisson structure  has a $sl(2)$ singularity at the origin if  
$j^2_0\hat f$ is R-equivalent to $x^2 - y^2$ or to $-x^2 - y^2$ 
(then $j^1_0\hat g = 0$ by (1.4)).  

\vskip .2cm

4) A Poisson structure  has an $A^+$ (resp. $A^-$) 
singularity at the origin if  
$j^2_0\hat f$ is R-equivalent to $x^2$ (resp. $-x^2$) 
and $j^1_0\hat g = 0$. An $A$ singularity is either $A^+$ or 
$A^-$ singularity.

\vskip .2cm

5) A Poisson structure  has a $N$ singularity at the origin if  
$j^2_0\hat f=0$ and $j^1_0\hat g = 0$. Within $N$ singularities 
we will study only $N^+$ and $N^-$ singularities - the cases  
where $j_0^2\hat g$ is $R$-equivalent to $x^2 + y^2$ and 
$x^2 - y^2$ respectively. 

\vskip .2cm

It is clear that any Poisson structure $P$ such that $P(0)=0, j^1_0P\ne 0$ 
has at the origin one of these 5 types of singularities. 
The singularity classes $V, so(3), sl(2), A$ and $N$ 
are related to the classification of 
3-dimensional Lie algebras since each of the singularity classes 
 is distinguished by a 
condition on the 1-jet at the origin of a Poisson structure, and the  
1-jet of any Poisson structure $P, P(0)=0,$ 
can be identified with a Lie algebra. It is not hard to check the 
following facts.

\vskip .2cm

a) The $so(3)$ and $sl(2)$ singularities correspond to the Lie algebras   
$so(3)$ and $sl(2)$ respectively (up to isomorphism). 

\vskip .2cm

b) The Lie algebras 
corresponding to the other singularity classes are isomorphic to a Lie 
algebra of the form $[e_1,e_2]=0, \ \  [e_1,e_3]=b_{1,1}e_1 + 
b_{1,2}e_2, \ \ \ [e_2,e_3]=b_{2,1}e_1 + 
b_{2,2}e_2$ with real parameters $b_{i,j}$.

\vskip .2cm

c) Let $B = (b_{i,j})$. 
The $V$ singularities are distinguished by the condition $trace B\ne 0$, 
the $A$ singularities by the condition $trace B = 0, \ det B\ne 0$ 
($det B > 0$ in the case of $A^+$ singularities and $det B < 0$ 
in the case of $A^-$ singularities), and 
the $N$ singularities  by the condition $trace B = 0, \ det B = 0, \ 
B\ne 0$. The singularity classes $N^+$ and $N^-$ cannot be 
distinguished in terms of the 1-jet of a Poisson structure (they 
are distinguished by a condition on $j_0^2P$).

\vskip .3cm 

The singularities $so(3)$ and $sl(2)$ are well known due to the works 
[Co1,Co2, We1].  
A Poisson structure having a $so(3)$ or $sl(2)$ singularity at the 
origin is formally equivalent to the linear Poisson structure 

$$\{x,y\} = z,  \ \{y,z\} = x, \ \{z,x\} = \pm y,$$ 

\noindent where the sign $+$ (resp. $-$) corresponds to the $so(3)$ 
(resp. $sl(2)$) singularity. In the case of $so(3)$ singularities 
this normal form also holds in the $C^\infty $ category.

The $so(3)$ and $sl(2)$ singularities are isolated and irremovable: 
if there is a 
family $P_\epsilon $, $\epsilon \in R^l$, of Poisson structures such 
that $P_0$ has a  
$so(3)$ (resp. $sl(2)$) singularity at the origin then there is a neighbourhood 
$W$ of $0\in R^l$ and a neighbourhood $U$ of $0\in R^3$ such that for any 
$\epsilon \in W$ the Poisson structure $P_\epsilon $ has a unique   
singular point in $U$ at which a $so(3)$ (resp. $sl(2)$) singularity holds.

\vskip .3cm 

The beginning of the classification of $V$ singularities 
can be found in the work [Du]. The classification 
reduces to the orbital classification of vector fields on a plane due to 
the following result: a local family $P_\epsilon $ of 
Poisson structures such that $P_0$ has a $V$ singularity at the 
origin is $C^\infty $ equivalent to a family of the form (1.5). 
In view of the correspondence between 3-dimensional Poisson structures 
and integrable 1-forms, this result is an analogous of the well-known 
Kupka phenomenon [Ku]: the local classification of integrable 1-forms 
$\omega $ on $R^3$ such that $d\omega (0)\ne 0$ reduces to the classification of 
arbitrary 1-forms on a plane.  

It follows 
that the $V$ singularities are also irremovable, but their geometry 
essentially differs from that for the $so(3)$ and $sl(2)$ singularities: 
if a Poisson structure has a $V$ singularity at the origin then it also has 
a $V$ singularity at each point of a smooth curve passing through the 
origin. A detailed study of the $V$ singularities is contained in 
section 2.

\vskip .3cm

The $A$ singularities are studied in section 3. These singularities 
can be removed by a small perturbation of an individual Poisson structure, 
but they are irremovable in 1-parameter families of Poisson structures. 
We give a normal form for any deformation of an algebraically isolated 
$A$ singularity. In particular, an individual 
Poisson structure having an algebraically isolated $A$  
singularity at the origin is formally equivalent to a Poisson structure 
of the form 
$$H(x,y,z)(z{\partial \over \partial x}\wedge {\partial \over \partial y} 
\pm  x{\partial \over \partial y}\wedge {\partial \over \partial z}  
\pm y^m{\partial \over \partial z}\wedge {\partial \over \partial x}), 
$$   
where $m\ge 2$ and $H$ is a nonvanishing function. This statement, 
expressed in 
terms of integrable 1-forms, can also be obtained as a 
corollary of the results of 
Moussu [Mo] or Malgrange [Ma].  
 
The geometry of a 
generic 1-parameter perturbation $P_\epsilon $ of a Poisson structure $P$ 
having a generic $A$ singularity at the origin is as follows. If 
$\epsilon < 0$ then $P_\epsilon $ contains no singular points near the 
origin; the Poisson structure $P=P_0$ has an isolated singular point which 
decomposes into two singular points near the origin as $\epsilon >0$. These  
singular points are both $sl(2)$ singularities of $P_\epsilon $  
if $P$ has an $A^-$ singularity at the origin. If $P$ has an 
$A^+$ singularity at the origin then    
one of the singular points is a $so(3)$ singularity and the other is a  
$sl(2)$ singularity. It is remarkable that no perturbation of arbitrary 
algebraically isolated $A$ singularity leads to a $V$ singularity whereas 
the $A$ singularities of the Lie algebras are adjoint to the 
$V$ singularities.

\vskip .3cm

The most difficult are the $N^+$ and $N^-$ 
singularities  studied in section 4. 
We prove that generic $N^+$ and $N^-$ singularities are irremovable in 
1-parameter families of Poisson 
structures. This is a bit surprising since, on the level of the Lie 
algebras, the $N$ singularities are typical only in 2-parameter families.  
We prove that there are three types of bifurcations in generic 
1-parameter families $P_\epsilon $ such that $P_0$ has a generic 
$N^+$ or $N^-$ singularity at the origin:

\vskip 0.2cm 

a) If $\epsilon \le 0$ then $P_\epsilon $  
has an isolated singular point which is a $so(3)$ singularity if 
$\epsilon <0$ and $N^+$ singularity if $\epsilon =0$. If $\epsilon 
>0$ then the set of singular points of $P_\epsilon $ 
consists of an isolated singular point and a closed curve. The isolated 
singular point is a $sl(2)$ singularity, and each point of the curve 
is a $V$ singularity. 

\vskip .2cm   

b) If $\epsilon \le 0$ then $P_\epsilon $  
has an isolated singular point which is a $sl(2)$ singularity if 
$\epsilon <0$ and $N^+$ singularity if $\epsilon =0$. If $\epsilon 
>0$ then the set of singular points of $P_\epsilon $ 
consists of an isolated singular point and a closed curve. The isolated 
singular point is a $so(3)$ singularity, and each point of the curve 
is a $V$ singularity. 

\vskip .2cm

c) The set of singular points of $P_\epsilon $ has the form (in suitable 
coordinates) $\{z=0, x^2 - y^2 - \epsilon = 0\}\cup \{x=y=z=0\}$. 
Any singular point except the origin is a $V$ singularity. The origin is 
a $sl(2)$ singularity if $\epsilon \ne 0$ and a $N^-$ singularity 
if $\epsilon = 0$. 

\vskip .2cm 

The type of the bifurcation and the type of appearing $V$ singularities 
(node, saddle, focus) can be determined in terms of the 3-jet of $P_0$.

\vskip .3cm 

\centerline {\bf 2. $V$-singularities.} 

\v 

Using the {\it curl} of a Poisson structure we  show in section 2.1 that 
the classification of $V$ singularities 
reduces to the orbital classification of vector fields on a plane (the 
analogous of the Kupka phenomenon [Ku]).  
Applying known results on the latter classification we obtain, in sections 
2.2, corollaries on normal forms, geometry and  
 bifurcations.  
 
 The $V$ singularities are never isolated - they hold 
at points of smooth curves. We distinguish hyperbolic 
$V$ singularities (node, saddle and focus)  and 
saddle-node $V$ singularities. The hyperbolic $V$ singularities  
are irremovable, and generic saddle-node $V$-singularities  
are irremovable in 1-parameter families (in generic 1-parameter 
families the saddle-node bifurcation holds). 
The results of section 2 continue the results of the paper 
[Du], where the nonresonant $V$ singularities were classified.

\v 

\centerline {\bf 2.1. Reduction to vector fields.}

\v 

 A volume form  $\Omega$ on $R^3$ induces  the isomorphism 
$b: Z\mapsto i_Z\Omega$ between $p$-vectors and $(3-p)$-forms on $R^3.$
The curl of a Poisson structure $P$ with respect to $\Omega $ 
is defined to be the vector field 
$X=b^{-1}(d(b(P))).$ It is known (see [DH], [We2]) 
that for any $\Omega $ the curl $X$ is an
infinitesimal symmetry of $P$, i.e., $[X, P]= 0$. If $P(0)=0$ then 
the vector $X(0)$  does not depend on $\Omega $. Computing the 
curl of the Poisson structure (1.3) with respect to the volume form 
$dx\wedge dy\wedge dz$, we obtain $X = -{\partial g\over \partial y}
{\partial \over \partial x} + {\partial g\over \partial x}
{\partial \over \partial y}$. This relation proves 
the following statement.  

\v

{\bf Proposition 2.1.} {\it For a Poisson structure $P$, $P(0)=0$, 
the origin is a $V$ singularity 
if and only if the curl of $P$ with respect to some 
(and then any)  volume form does not vanish at the origin.} 

\v

Assume that we have a family $P_\epsilon $ of Poisson structures such 
that $P_0$ has a $V$ singularity at the origin. Let 
$X_\epsilon $ be the curl of $P_\epsilon $ with respect to a fixed 
volume form $\Omega $.  
By Proposition 2.1 
there exists a coordinate system (depending smoothly on $\epsilon $) 
such that $X_\epsilon = {\partial \over \partial y}$. In the 
3-dimensional case the Jacobi formula implies $X_\epsilon \wedge P_\epsilon = 0.$ 
It follows from this relation and the relation 
$[X_\epsilon , P_\epsilon ]= 0$ that in the choosen coordinate system 
$P_\epsilon $ has the form (1.5). Returning to the characterization 
of $V$ singularities in terms of the 1-jet of a Poisson structure, 
we conclude that the sum of the eigenvalues of the vector field $v_0$ in 
(1.5) is different from zero.      
So, we have proved the following 
statement.   

\v 

{\bf Proposition 2.2.} {\it Any local family $P_\epsilon $ of  
Poisson structures on $R^3$ 
such that $P_0$ has a $V$-singularity 
at the origin is  equivalent to a family of the form   
$${\partial \over \partial y}\wedge v_\epsilon , \ \ v_\epsilon  = \alpha _\epsilon (x,z){\partial \over  
\partial x} + \beta _\epsilon (x,z){\partial \over \partial z},$$ 
where the sum of the eigenvalues of the linearization at the origin 
of the vector field $v_0$ is different from zero.}

\v

Consider now two families of Poisson structures 
$P_\epsilon = {\partial \over \partial y}\wedge v_\epsilon $ and   
$\tilde P_\epsilon = {\partial \over \partial y}\wedge \tilde v_\epsilon $, 
where  
$v_\epsilon  = \alpha _\epsilon (x,z){\partial \over 
\partial x} + \beta _\epsilon (x,z){\partial \over \partial z}$ and 
$\tilde v_\epsilon  = \tilde \alpha _\epsilon (x,z){\partial \over 
\partial x} + \tilde \beta _\epsilon (x,z){\partial \over \partial z}$. 
Assume that the families $v_\epsilon $ and $\tilde v_\epsilon $ are 
orbitally equivalent, i.e., there exists a family $\phi _\epsilon $ of 
local diffeomorphisms of the $(x,z)$ plane such that 
$(\phi _\epsilon )_*v_\epsilon = 
h_\epsilon \tilde v_\epsilon $, where $h_\epsilon $ is a family of 
functions such that $h_0(0)\ne 0$. Then the diffeomorphism 
$(x,z)\mapsto \phi _\epsilon (x,z), y\mapsto y/h_\epsilon $ brings 
$P_\epsilon $ to $\tilde P_\epsilon $. So, the orbital equivalence of 
the families  $v_\epsilon $ and $\tilde v_\epsilon $ implies the 
equivalence of the families $P_\epsilon $ and $\tilde P_\epsilon $. 
The inverse statement holds under the assumption that the fields 
$v_0$ and $\tilde v_0$ have isolated singularities at the origin. 
In fact, let $\psi _\epsilon $ be a family of diffeomorphisms sending 
$P_\epsilon $ to $\tilde P_\epsilon $. Since the $y$ axis is the set of 
singular points for both $P_0$ and $\tilde P_0$, the 
diffeomorphism $\psi _0$ sends the plane $y=0$ to a surface 
transversal to the $y$ axis. This property remains true for the 
diffeomorphism $\psi _\epsilon $, when $\epsilon $ is small. 
Therefore there exists a family of functions 
$g_\epsilon (x,z)$ such that the superposition $\mu _\epsilon $ 
of $\psi _\epsilon $ 
with the diffeomorphism  
$\nu _\epsilon  : (x,y,z)\mapsto (x, y - g_\epsilon (x,z), z)$ preserves the 
plane $y=0$ for each small enough $\epsilon $. The diffeomorphism $\nu _\epsilon $ 
is a symmetry of the Poisson structures $P_\epsilon $ and 
$\tilde P_\epsilon $. It follows that the superposition $\mu _\epsilon $ 
also brings $P_\epsilon $ to $\tilde P_\epsilon $. It is easy to check that 
the 
restriction of $\mu _\epsilon $ to the plane $y=0$ brings  
$v_\epsilon $ to $\tilde v_\epsilon $ multiplied by a nonvanishing 
function, i.e., the families $v_\epsilon $ and $\tilde v_\epsilon $ 
are orbitally equivalent. We have proved the following statement reducing 
the classification of $V$ singularities to the orbital classification of 
vector fields. 

\v  

{\bf Proposition 2.3.} {\it  Let 
$v_\epsilon = \alpha _\epsilon (x,z){\partial \over 
\partial x} + \beta _\epsilon (x,y){\partial \over \partial z}$ 
and 
$\tilde v_\epsilon = \tilde \alpha _\epsilon (x,z){\partial \over 
\partial x} + \tilde \beta _\epsilon (x,y){\partial \over \partial z}$ 
be two local families of vector fields such that the vector fields 
$v_0$ and $\tilde v_0$ have isolated singular point at the origin. 
The family of Poisson structures $P_\epsilon = 
{\partial \over \partial y}\wedge v_\epsilon $ is equivalent to 
the family $\tilde P_\epsilon = 
{\partial \over \partial y}\wedge \tilde v_\epsilon $ if and only if 
the family $v_\epsilon $ is orbitally equivalent to the family 
$\tilde v_\epsilon $. } 

\v 

\centerline {\bf 2.2. Normal forms and geometry of $V$ singularities.} 

\v 
The results of this section are direct corollaries of Propositions 
2.2 and 2.3 and results on the orbital classification of 
vector fields on a plane, see [AI, Bo]. At first we distinguish node, saddle, focus and 
saddle-node $V$ singularities. 
Assume that a Poisson structure $P$ has a $V$ singularity at the origin. 
By Proposition 2.2, $P$ is equivalent to a Poisson structure of the form 
${\partial \over \partial y}\wedge v$, where $v$ is a vector field on 
the $(x,z)$ plane, $v(0)=0$.  
Let $\lambda _1$ and $\lambda _2$ be the eigenvalues of the 
linearization of the vector field $v$ in this normal form. By Proposition 
2.3 they are invariantly related to $P$ up to multiplication by 
a common factor, and we will say that $\lambda _1$ and $\lambda _2$ are 
the eigenvalues of $P$.  Note that $\lambda _1 + \lambda _2 \ne 0$, see Proposition 
2.2.

We will say that a $V$ singularity is hyperbolic if $\lambda _1\ne 0$ and 
$\lambda _2\ne 0$. Since $\lambda _1 + \lambda _2\ne 0$ this means that the 
spectrum of the linearization of the vector field $v$ does not 
intersect the imaginary axis.

Within hyperbolic $V$ singularities we distinguish node $V$ singularities 
corresponding to the case where $\lambda _1$ and 
$\lambda _2$ are real nonzero numbers of the same sign, 
saddle $V$ singularities corresponding to the case where 
$\lambda _1$ and 
$\lambda _2$ are real nonzero numbers of different signs, and 
focus $V$ singularities corresponding to the case where $\lambda _{1,2} = 
a\pm ib$, $a\ne 0$, $b\ne 0$.

If a $V$ singularity is not hyperbolic then one of the eigenvalues 
is equal to zero and the second is different from zero.    
In this case we will say that the origin is a saddle-node $V$ singularity. 
Within saddle-node $V$ singularities there exists a degeneration of 
infinite codimension corresponding to the case where the origin is not 
algebraically isolated singular point of the vector field $v$  
(this is impossible for hyperbolic $V$ singularities). 
If such a degeneration holds we will say that the origin is exclusive 
saddle-node $V$-singularity.  

\v 

{\bf Theorem 2.4.} {\it Assume that $P$ has a hyperbolic or nonexclusive 
saddle-node $V$-singularity at the origin. Then 
the set of singular points of $P$ in a small enough 
neighbourhood of the origin is a smooth curve $\gamma $. 
The germs of $P$ at points of $\gamma $ are $C^\infty $ 
equivalent to the germ of $P$ at the origin and $C^\infty $ 
equivalent to one of the following normal forms:} 

$$\{x,z\}=0, \ \{x,y\} = z, \ \{z,y\} = \theta x + z,\ \ \theta \in R-\{0\};  
\eqno (2.1)$$

$$\{x,z\}=0, \ \{x,y\} = Nx + \delta z^N, \ \{z,y\} = z, \ 
N\in \{1, 2, 3, 4...\}, \ \delta \in \{0,1\};\eqno (2.2)$$

$$\{x,z\}=0, \ \{x,y\} = -{p\over q}x + \delta x^{q+1}z^p + \delta ax^{2q+1}z^{2p}, 
\ \{z,y\} = z, \ \delta \in \{1,-1,0\}, a\in R;   
\eqno (2.3)$$

$$\{x,z\}=0, \ \{x,y\} = \delta ^{p+1}x^{p+1} + ax^{2p+1}, \ \{z,y\} = z, \ 
p\ge 1, \ \delta \in \{1,-1\}, \ a\in R. \eqno (2.4)$$

\v 

The normal form (2.1) holds for focus $V$ singularities and nonresonant 
node or saddle $V$ singularities, i.e., for node $V$ singularities such 
that neither $\lambda _1/\lambda _2$ nor 
$\lambda _2/\lambda _1$ is an integer  number and for 
saddle $V$ singularities such that $\lambda _1/\lambda _2$ is 
not a rational number. Here $\lambda _1$ and $\lambda _2$ are the 
eigenvalues of $P$. The normal form (2.2) holds for resonant 
node $V$ singularities, the normal form (2.3) holds for resonant 
saddle $V$ singularities (in this normal form $p$ and $q$ are positive 
integer numbers and $p/q$ is an irreducible fraction), and the 
normal form (2.4) holds for nonexclusive saddle-node $V$ singularities.  
  
Since the hyperbolic singular points of vector fields are irremovable, 
the hyperbolic $V$ singularities are irremovable under a small perturbation 
of a Poisson structure. Namely,    
if $P$ has a node (resp. saddle, focus) $V$ singularity at the origin  
and $P_\epsilon $ is a family of 
Poisson structures, $\epsilon \in R^l$, such that $P_0=P$ then there exist 
a neighbourhood of the origin $U\subset R^3$ 
and a neighbourhood of the origin $W\subset R^l$ such that for any 
$\epsilon \in W$ the set of singular points of $P_\epsilon $ in $U$ 
is a smooth curve $\gamma _\epsilon $. The family $\gamma _\epsilon $ 
depends smoothly on $\epsilon $ and each point of $\gamma _\epsilon $ is 
also a node (resp. saddle, focus) $V$ singularity of $P_\epsilon $.

The saddle-node $V$ singularities are irremovable only in 1-parameter 
families of Poisson structures. Using the correspondence between $V$ 
singularities and vector fields on a plane, we will say that a 
saddle-node $V$ singularity is generic if $p=1$ in the normal form 
(2.4). If $P$ has a generic saddle-node singularity at the origin and 
$P_\epsilon $ is a 1-parameter deformation of $P$, then in suitable 
coordinate system (depending smoothly on $\epsilon $) the 
2-jet of the family $P_\epsilon $ has the form 
$\{x,z\}=0, \ \{x,y\} = f_0(\epsilon ) + f_1(\epsilon )x + x^2, \  
\{z,y\} = z$. We will say that $P_\epsilon $ is a generic deformation 
of $P$ if $f_0^\prime (0)\ne 0$. In this case the following (up to the 
change $\epsilon \to -\epsilon $) saddle-node bifurcation holds (the 
analogous of the well known saddle-node bifurcation for vector fields).  
There exists a neighbourhood of the origin $U\subset R^3$ 
and a neighbourhood of the origin $W\subset R$ such that    
if $\epsilon \in W$ is a negative number then $U$ contains no 
singular points of $P_\epsilon $ and if $\epsilon \in W$ is a 
positive number then the set of singular points of $P_\epsilon $ in 
$U$ consists of two 
disjoint smooth curves.  The points of the first  curve are saddle 
$V$ singularities and the points of the second curve are node  
$V$ singularities of $P_\epsilon $.   

\v 

\centerline {\bf 3. $A$ singularities.}

\v 

This section contains results on algebraically isolated $A$ singularities. 
Recall that a Poisson structure 
$\{x ,y \}=B(x,y,z) , \ \{y ,z \} = C(x,y,z),\ \{z ,x \} = D(x,y,z)$ has 
an algebraically isolated singularity at the origin if the factor ring 
of the ring of all formal series over the ideal generated by the 
formal series of the functions $B,C$ and $D$ has finite dimension over $R$. 
It is clear that an $A$ singularity is algebraically isolated if and only if 
the formal series $\hat f_0(x,y)$ in the normal form (1.3)    
is $R$-equivalent to $\pm x^2\pm y^{m+1}$ for some $m\ge 2$. 
It follows that all $A$ singularities, except degenerations of 
infinite codimension, are algebraically isolated. We will say that an 
$A$ singularity is generic if $m=2$.     

\v 

{\bf Theorem 3.1.} {\it Let $P_\epsilon $ be a local family of Poisson 
structures of the form (1.2) such that 
$P_0$ has an algebraically isolated $A$ singularity at the origin 
$p_0$. There is a family of smooth changes of 
coordinates parametrized by $\epsilon $ such that in the new coordinate system $P_\epsilon $ 
has the form   
$$H_\epsilon (x,y,z)(z{\partial \over \partial x}\wedge 
{\partial \over \partial y} \pm x{\partial \over \partial y}\wedge 
{\partial \over \partial z} + (\delta y^m + 
\sum _{i=0}^{m-2}h_i(\epsilon )y^i){\partial \over \partial z}\wedge 
{\partial \over \partial x})\eqno (3.1)$$ 
modulo flat field of 2-vectors at the point $p_0$, 
where $m\ge 2$, $H_0(0) > 0$, $h_0(0) = \dots =  
h_{m-2}(0) = 0$,     
$\delta = \pm 1$ if $m$ is odd and  $\delta = 1$ if 
$m$ is even.}

\v 

The sign $+$ (resp. $-$) in (3.1) corresponds to the case of $A^+$ 
(resp. $A^-$) singularities. Note that in general the point $p_0$ 
is the origin of the coordinate system of the normal form (3.1)  
only if $\epsilon =0$. 
 
In particular, any 1-parameter unfolding of a generic $A$ singularity reduces to 
the formal normal form 
$$H_\epsilon (x,y,z)(z{\partial \over \partial x}\wedge 
{\partial \over \partial y} \pm x{\partial \over \partial y}\wedge 
{\partial \over \partial z} + (y^2 + h_0(\epsilon ))
{\partial \over \partial z}\wedge 
{\partial \over \partial x}).\eqno (3.2)$$ 

A 1-parameter unfolding of a generic  
$A$ singularity will be called generic 
if the function $h_0$ in this normal form satisfies the condition 
$h_0^\prime (0)\ne 0$. The following theorem says that $A$ singularities 
are irremovable in 1-parameter families of Poisson structures and a 
generic individual $A$ singularity decomposes, under a small perturbation, 
into two singular points; the type of singularity at these points 
depends on the type of singularity of $P$ at the origin ($A^+$ or $A^-$).

\v 

{\bf Theorem 3.2.} {\it Let $P_\epsilon $ be a generic 1-parameter 
unfolding of a generic Poisson structure $P$ having an $A$ singularity 
at the origin. Then the following bifurcation (up to the change 
$\epsilon \to -\epsilon $) holds in a neighbourhood 
of the origin $U\subset R^3$ 
and in a neighbourhood of the origin $W\subset R$. If 
$\epsilon \in W$ and $\epsilon < 0$ then $U$ contains 
no singular points of $P_\epsilon $, if $\epsilon \in W$ and $\epsilon > 0$ 
then $U$ contains two singular points of $P_\epsilon $. 
If the singularity of $P$ has the type 
$A^+$ then one of these singular points is a $so(3)$ singularity and the 
other is a $sl(2)$ singularirty. If the singularity of $P$ has the type  
$A^-$ then the two singular points are $sl(2)$ singularities.} 

\v 

Theorem 3.2 is an easy corollary of the normal form (3.2). 
One also can obtain, using the normal form (3.1), that any (not 
necessarily generic) algebraically isolated $A^-$ singularity 
decomposes, under a small perturbation, into $sl(2)$ singularities, and 
any algebraically isolated $A^+$ singularity decomposes into 
$so(3)$ and $sl(2)$ singularities. It follows that there is no 
adjacency between algebraically isolated $A$ singularities and $V$ 
singularities. Note that such an adjacency holds for Lie algebras - 
a suitable perturbation of any Lie algebra of the type $A$   
(in the class of Lie algebras) leads to a Lie algebra of the type $V$. 

\v 

{\bf Proof of Theorem 3.1.} We will use Theorem 1.1, i.e., 
the formal normal form (1.3)-(1.4).  Since 
$P_0$ has an algebraically isolated singularity at the origin then 
the 1-form $d\hat f_0(x,y)$ also has an algebraically isolated singularity 
at the origin of the $(x,y)$ plane. Using this observation and results 
of the paper [Mo] we conclude that (1.4) implies that 
$\hat g_\epsilon = \hat \lambda _\epsilon \circ \hat f_\epsilon $, where 
$\hat \lambda _\epsilon $ is a family of formal series in one variable 
depending smoothly on $\epsilon $. 

Now we will show that $P_\epsilon $ admits a family of Casimir 
functions, i.e., there exists a family $\hat C_\epsilon = 
\hat C_\epsilon (x,y,z)$ of nonzero formal series such that 
$\hat C_\epsilon (0) = 0$ and $\{\hat C_\epsilon , x\} = 
\{\hat C_\epsilon , y\} = \{\hat C_\epsilon , z\} = 0$. The latter 
relation means that $P_\epsilon (d\hat C_\epsilon , d\hat h) = 0$ for 
any formal series $\hat h$. 

Taking into account the relation 
$\hat g_\epsilon = \hat \lambda _\epsilon \circ \hat f_\epsilon $ it is 
natural to seek for a family of Casimir functions in the form 
$\hat C_\epsilon (x,y,z) = \hat G_\epsilon (z, \hat f_\epsilon )$, where 
$\hat G_\epsilon = \hat G_\epsilon (z,w)$ is a family of formal 
series in two variables. A simple computation shows that 
$$P_\epsilon (d\hat C_\epsilon , dz) = 0, \ 
P_\epsilon (d\hat C_\epsilon , dx) = \hat Q_\epsilon {\partial \hat f_\epsilon \over 
\partial y}, \ P_\epsilon (d\hat C_\epsilon , dy) = 
-\hat Q_\epsilon {\partial \hat f_\epsilon \over 
\partial x},$$ where 
$$\hat Q_\epsilon = (1 + z\hat \lambda _\epsilon ^\prime  (\hat f_\epsilon ))
{\partial \hat G_\epsilon \over \partial z}(z,\hat f_\epsilon ) - 
z{\partial \hat G_\epsilon \over \partial w}(z, \hat f_\epsilon ).$$ 
The equation $\hat Q_\epsilon = 0$ has a formal solution $\hat G_\epsilon 
(z,w)$ of the form 
$$\hat G_\epsilon (z,w) = w + z^2/2 + a(\epsilon )w^2 + \hat R_\epsilon (z,w),$$ 
where $\hat R_\epsilon $ is a formal series with zero 2-jet for each fixed 
$\epsilon $.  Therefore $P_\epsilon $ 
admits a family of Casimir functions $\hat C_\epsilon (x,y,z)$ such that 
$$\hat C_0(x,y,z) = \hat f_0(x,y) + z^2/2 + a(0)\hat f_0^2(x,y) +  
R_0(z,\hat f_0(x,y)). $$ 
The series $\hat f_0(x,y)$ is $R$-equivalent to $\pm {x^2\over 2} + 
\delta {y^{m+1}\over m+1}$, where $m\ge 2$, 
$\delta = 1$ if $m$ is even and  $\delta 
\in \{-1,1\}$ if $m$ is odd. It follows that the series $\hat C_0(x,y,z)$ 
is $R$-equivalent to ${z^2\over 2} \pm {x^2\over 2} + 
\delta {y^{m+1}\over m+1}$. Let $C_\epsilon (x,y,z)$ be 
a family of smooth functions such that the formal series of 
$C_\epsilon (x,y,z)$ is equal to $\hat C_\epsilon (x,y,z)$. 
Then there exists a new smooth system  
of coordinates $\tilde x = \tilde x_\epsilon , \tilde y = \tilde y_\epsilon , 
\tilde z = \tilde z_\epsilon $  such that $$C_\epsilon (x,y,z) = 
\tilde C_\epsilon (\tilde x, \tilde y, \tilde z) =   
\alpha (\epsilon ) + {\tilde z^2\over 2} \pm {\tilde x^2\over 2} + 
\delta {\tilde y^{m+1}\over m+1} + \sum _{i=1}^{m-1}h_i(\epsilon )
\tilde y^i,$$ 
where $h_i(0)=0$, see [AVG]. Now Theorem 3.1 easily follows from  
the observation that the functions 
$$\{\tilde C_\epsilon (\tilde x, \tilde y, \tilde z), \tilde x\}, 
\{\tilde C_\epsilon (\tilde x, \tilde y, \tilde z), \tilde y\}, 
\{\tilde C_\epsilon (\tilde x, \tilde y, \tilde z), \tilde z\}$$ are flat 
at the point $x=y=z=0$.

\v

\centerline {\bf 4. $N^+$ and $N^-$ singularities.}

\v 

In this section we show that the $N^+$ and $N^-$ singularities are 
typical in 1-parameter families of Poisson structures (though they 
are typical only in 2-parameter families of Lie algebras) and 
analyse typical bifurcations in 1-parameter families $P_\epsilon $ 
such that $P_0$ has an $N^+$ or $N^-$ singularity at the origin. 
The analysis of bifurcations is based on the following proposition. 

\v 

{\bf Proposition 4.1.} {\it Any local family of Poisson structures $P_\epsilon $ 
such that $P_0$ has a $N^+$ or $N^-$ singularity at the origin is  
equivalent to a family of the form 
$$\{x,y\} = z, \ \{y,z\} = A_\epsilon (x,y)(z + C_\epsilon (x,y)) + 
z^2Q_{1,\epsilon }(x,y,z),$$ 
$$\{z,x\} = B_\epsilon (x,y)(z + C_\epsilon (x,y)) + 
z^2Q_{2,\epsilon }(x,y,z), \eqno (4.1)$$ 
where the formal series of the functions $A_\epsilon $, $B_\epsilon $ 
and $C_\epsilon $ have the form 
$$\hat A_\epsilon (x,y) = x\left (\lambda _0(\epsilon ) + \lambda _1
(\epsilon )(x^2\pm y^2) + \lambda _2(\epsilon )(x^2\pm y^2)^2 + \cdots 
\right ),$$ 
$$\hat B_\epsilon (x,y) = \pm y\left (\lambda _0(\epsilon ) + \lambda _1
(\epsilon )(x^2\pm y^2) + \lambda _2(\epsilon )(x^2\pm y^2)^2 + \cdots 
\right ),\eqno (4.2)$$  
$$\hat C_\epsilon (x,y) = \mu _0(\epsilon ) + \mu _1
(\epsilon )(x^2\pm y^2) + \mu  _2(\epsilon )(x^2\pm y^2)^2 + \cdots ,   
$$ 
and the functions $\lambda _0(\epsilon )$ and  $\mu  _0(\epsilon )$  
satisfy the conditions 
$$\lambda _0(0)\ne 0, \ \mu _0(0) = 0. \eqno (4.3)$$ 
The signs $+$ (resp. $-$) in (4.2) correspond to the case where $P_0$ 
has a $N^+$ (resp. $N^-$) singularity at the origin. }

\v

This proposition will be proved at the end of the section. 
In terms of the normal form (4.1)-(4.2) we define generic $N^+$ or 
$N^-$ singularity and its generic 1-parameter unfolding. A $N^+$ or 
$N^-$ singularity will be called generic if $\mu _1(0)\ne 0$, and 
a 1-parameter unfolding will be called generic if $\mu _0^\prime (0)
\mu _1(0)\ne 0$. 

The bifurcations in generic 1-parameter unfoldings of generic 
$N^+$ and $N^-$ singularities depend on the signs of the numbers 
$$\kappa _1 = \lambda _0(0)\mu _1(0), \ \kappa _2 = \lambda _0^2(0) - 
8\kappa _1.$$ The main result of this section is the following theorem. 

\v 

{\bf Theorem 4.2.} {\it Let $P$ be a Poisson structure with a generic 
$N^+$ or $N^-$ singularity at the origin, and let $P_\epsilon $ be 
a generic 1-parameter unfolding of $P$. There exists a coordinate 
system $\tilde x = \tilde x_\epsilon , \tilde y = \tilde y_\epsilon  , 
\tilde z = \tilde z_\epsilon $, $\tilde \epsilon = \epsilon Q(\epsilon )$, 
$Q(0)\ne 0$, such that the following statements hold (locally near the 
point $\tilde x = \tilde y = \tilde z = \tilde \epsilon = 0$).  

1) The set of singular points of $P_\epsilon $ has the form 
$$\{\tilde z = 0, \ \tilde x^2 + \tilde y^2 - \tilde \epsilon = 0\}\cup 
\{\tilde x = \tilde y = \tilde z = 0\}$$ in the case of $N^+$ 
singularity and   
$$\{\tilde z = 0, \ \tilde x^2 - \tilde y^2 - \tilde \epsilon = 0\}\cup 
\{\tilde x = \tilde y = \tilde z = 0\}$$ in the case of $N^-$ 
singularity. 

2) The Poisson structure $P_\epsilon $ has a $V$ singularity at any 
singular point except the point $\tilde x = \tilde y = \tilde z = 0$.  
Any $V$ singularity is a node $V$ singularity if $\kappa _2\ge 0$ and 
$\kappa _1 >0$, a saddle $V$ singularity if $\kappa _2\ge 0$ and 
$\kappa _1 <0$, and a focus $V$ singularity if $\kappa _2<0$.

3) If $P$ 
has a $N^+$ singularity at the origin and $\kappa _1 >0$ (resp. 
$\kappa _1<0$) then the point $\tilde x = \tilde y = \tilde z = 0$ is 
a $so(3)$ (resp. $sl(2)$) singularity of $P_\epsilon $ if $\epsilon >0$  
and a $sl(2)$ (resp. $so(3)$) singularity if $\epsilon <0$.

4) If $P$ has a $N^-$ singularity at the origin then for any 
$\epsilon \ne 0$ the Poisson structure $P_\epsilon $ has a $sl(2)$ 
singularity at the point $\tilde x = \tilde y = \tilde z = 0$.} 

\v 

This theorem implies that there are three types of typical bifurcations: 

a) $\epsilon <0$: a unique singular point - a $so(3)$ singularity; 
$\epsilon = 0$:  a unique singular point - a $N^+$ singularity such 
that $\kappa _1 >0$; $\epsilon  >0$: a $sl(2)$ singularity at an 
isolated singular point and a closed curve of $V$ singularities which 
are either nodes (if $\kappa _2\ge 0$) or focuses (if $\kappa _2 <0$). 

b) $\epsilon <0$: a unique singular point - a $sl(2)$ singularity; 
$\epsilon = 0$:  a unique singular point - a $N^+$ singularity such 
that $\kappa _1 <0$; $\epsilon  >0$: a $so(3)$ singularity at an 
isolated singular point and a closed curve of saddle $V$ singularities. 

c) The set of singular points consists of a point $p$  
and two curves given by the equation 
$x^2 - y^2 = \epsilon , z = 0$ in a coordinate system centered at the 
point $p$. The point $p$ is a $sl(2)$ singularity if 
$\epsilon \ne 0$ and an $N^-$ singularity if $\epsilon =0$. All other 
singular points are node or saddle or focus 
$V$ singularities dependently on the signs of the numbers $\kappa _1$ and 
$\kappa _2$ only.

\v 

{\bf Proof of Theorem 4.2.} We will restrict ourselves to the case of 
$N^+$ singularities such that $\kappa _1 > 0$ (the proof for 
the other cases is similar). Throughout the proof we use the normal 
form (4.1)-(4.2). We will assume that $\mu _0^\prime (0)
\mu _1(0)<0$ (for if not one can change $\epsilon $ by $\tilde \epsilon = 
-\epsilon $). The first statement of Theorem 4.2 obviously 
follows from the normal form, and the second and the third statements 
are corollaries of the following calculations.

a) Assume that $\epsilon < 0$. The equation $C_\epsilon (x,y) = 0$ has no solutions 
near the origin, therefore the point $x=y=z=0$ is the only singular point 
of $P_\epsilon $. The linearization of $P_\epsilon $ at this point has 
the form 
$$\{x,y\} = z, \ \{y,z\} = \mu _0(\epsilon )\lambda _0(\epsilon )x, \ 
\{z,x\} = \mu _0(\epsilon )\lambda _0(\epsilon )y.\eqno (4.4)$$ 
The relations $\epsilon <0, \mu _0^\prime (0)\mu _1(0) < 0$ and 
$\mu _1(0)\lambda _0(0) >0$ imply that $\mu _0(\epsilon )
\lambda _0(\epsilon ) > 0$, and it follows that $P_\epsilon $ 
has a $so(3)$ singularity at the point $x=y=z=0$. 

b) Assume that $\epsilon > 0$.  The set of singular points of 
$P_\epsilon $ consists of the point $x=y=z=0$ and the closed 
curve $\Gamma _\epsilon $ given by the equations $z=0, C_\epsilon 
(x,y) = 0$. The linearization of $P_\epsilon $ at the point 
$x=y=z=0$ has the form (4.4), but now $\epsilon >0$ and $\mu _0(\epsilon )
\lambda _0(\epsilon ) < 0$. Therefore $P_\epsilon $ has a $sl(2)$ 
singularity at the point $x=y=z=0$. 

These calculations prove the third statement of Theorem 4.2. 
To prove the second statement, 
take a point $p$ of the curve $\Gamma _\epsilon $. The linearization 
of $P_\epsilon $ at $p$ has the form     
$$\{x,y\} = z, \ \{y,z\} = A_\epsilon (p)\left ({\partial C_\epsilon 
\over \partial x}(p)x + {\partial C_\epsilon 
\over \partial y}(p)y + z\right ), $$  
$$\{z,x\} = B_\epsilon (p)\left ({\partial C_\epsilon 
\over \partial x}(p)x + {\partial C_\epsilon 
\over \partial y}(p)y + z\right ).\eqno (4.5)$$ 
The Jacobi identity implies 
$$ A_\epsilon (p){\partial C_\epsilon 
\over \partial y}(p) = B_\epsilon (p){\partial C_\epsilon 
\over \partial x}(p).\eqno (4.6)$$ 
Assume that $B_\epsilon (p)\ne 0$. Let $u = A_\epsilon (p)x + B_\epsilon 
(p)y$. Taking into account (4.6) we obtain that in the coordinate 
system $x,z,u$ the Poisson structure (4.5) takes the form 
$$\{u,z\} = 0, \ \{x,u\} = B_\epsilon (p)z, \ \{x,z\} = 
- B_\epsilon (p)z - {\partial C_\epsilon \over \partial y}(p)u. 
\eqno (4.7)$$ Now we see that the condition 
$B_\epsilon (p)\ne 0$ implies that 
$P_\epsilon $ has a $V$ singularity at $p$. On the other hand, the same 
conclusion holds if $A_\epsilon (p)\ne 0$. In this case we use the 
coordinate system $y,z,u$ in which (4.5) takes the form  
$$\{u,z\} = 0, \ \{y,u\} = -A_\epsilon (p)z, \ \{y,z\} = 
A_\epsilon (p)z + {\partial C_\epsilon \over \partial x}(p)u.  
$$
Note now that for any point $p\in \Gamma _\epsilon $ either 
$A_\epsilon (p)\ne 0$ or $B_\epsilon (p)\ne 0$. Therefore $P_\epsilon $ 
has a $V$ singularity at any point of the curve $\Gamma _\epsilon $.

The type of $V$ singularity (saddle, node, focus) is the same for 
all points of $\Gamma _\epsilon $ (see section 2), therefore for 
its determination it suffices to analyse the linear approximation 
at a single point of $\Gamma _\epsilon $, for example, at a point 
$q$ which is the intersection of $\Gamma _\epsilon $ with the 
semiaxis $\{x=z=0, y>0\}$. The point $q$ has the coordinates 
$(0, r\sqrt \epsilon + o(\sqrt \epsilon ), 0)$, where 
$$r = \sqrt {-{\mu _0^\prime (0)\over \mu _1(0)}}.$$ 
The function $B_\epsilon $ does not vanish at $q$, therefore the 
linear approximation of $P_\epsilon $ at $q$ is isomorphic to 
(4.7). Note that 
$$B_\epsilon (q) = r\lambda _0(0)\sqrt \epsilon + o(\sqrt \epsilon ), $$ 
$${\partial C_\epsilon \over \partial y}(q) = 2r\mu _1(0)\sqrt \epsilon  
+ o(\sqrt \epsilon ).$$ 
These relations allow us to write (4.7) in the form 
$r\sqrt \epsilon {\partial \over \partial x}\wedge v_0 + o(\sqrt \epsilon )$, 
where 
$$v_0 = \lambda _0(0)z{\partial \over \partial u} - 
(\lambda _0(0)z + 2\mu _1(0)u){\partial \over \partial z}.$$ 
It follows that the type of the $V$ singularity of $P_\epsilon $ at 
$q$ is determined by the type of singular point of the vector field 
$v_0$. Computing the eigenvalues of $v_0$ we obtain the second statement 
of Theorem 4.2.

\v 

{\bf Proof of Proposition 4.1.} The idea of the proof is as follows. At 
first we prove that in suitable smooth coordinate system the curl of 
$P_\epsilon $ vanishes at the origin for all $\epsilon $. After this 
we reduce $P_\epsilon $ to the formal normal form (1.3)-(1.4) and simplify 
this normal form using the condition on the curl. This simplification 
allows us to obtain the smooth normal form (4.1)-(4.2).

Take a coordinate system $x,y,z$ (depending 
on $\epsilon $) such that $P_\epsilon $ has the form 
$$\{x,y\} = z,\ \{y,z\} = c_\epsilon (x,y) + za_\epsilon (x,y) + 
z^2Q_{1,\epsilon }(x,y,z),$$ $$\{z,x\} = d_\epsilon (x,y) + zb_\epsilon (x,y) + 
z^2Q_{2,\epsilon }(x,y,z).\eqno (4.8)$$ 
The Jacobi identity restricted to the plane $z=0$ gives the relation 
$$c_\epsilon (x,y)b_\epsilon (x,y) =  d_\epsilon (x,y)a_\epsilon (x,y). 
\eqno (4.9)$$ It is easy to see that the condition that $P_0$ has an $N^+$ 
or $N^-$ singularity at the origin implies that the function 
$a_0(x,y)$ and $b_0(x,y)$ are differentially independent. Therefore 
the relation (4.9) implies the existence of a family $s_\epsilon (x,y)$ 
such that 
$$c_\epsilon (x,y) = s_\epsilon (x,y)a_\epsilon (x,y), \ 
d_\epsilon (x,y) = s_\epsilon (x,y)b_\epsilon (x,y).\eqno (4.10)$$ 
Consider the family $p_\epsilon $ of points with coordinates $x,y,z$ 
such that 
$$z=0, \ a_\epsilon (x,y) = 0, \ b_\epsilon (x,y) = 0.$$ It is clear that 
$p_\epsilon $ depends smoothly on $\epsilon $ and that $p_\epsilon $ is 
a singular point of $P_\epsilon $. The linearization of $P_\epsilon $ 
at the point $p_\epsilon $ has the form 
$$\{x,y\} = z, \ \{y,z\} = s_\epsilon (p_\epsilon )(e_1x + e_2y), \ 
 \{z,x\} = s_\epsilon (p_\epsilon )(e_3x + e_4y).\eqno (4.11)$$ 
It is easy to check that the Jacoby identity implies that any linear 
Poisson structure of the form (4.11) has the curl vanishing at the 
point $x=y=z=0$. Therefore the curl of $P_\epsilon $ vanishes at any 
singular point $p_\epsilon $. 

There is no loss of generality to assume that $p_\epsilon = 0$ for 
all $\epsilon $, i.e., that the curl of $P_\epsilon $ vanishes 
at the origin, and it what follows we assume that this condition holds. 
By Theorem 1.1 the family $P_\epsilon $ can be reduced, by a formal 
change of coordinates, to the form (1.3)-(1.4). It is easy to see that any 
formal change of coordinates of the form $$\phi _\epsilon : (x,y)\to 
(\phi _{1,\epsilon }(x,y), \phi _{2,\epsilon }(x,y))\eqno (4.12)$$ 
preserving the 
volume form $dx\wedge dy$ also preserves the normal form (1.3)-(1.4)  
up to the change $\hat f_\epsilon \to 
\hat f_\epsilon \circ \phi _\epsilon ,  
\hat g_\epsilon \to \hat g_\epsilon \circ \phi _\epsilon $. 
The fact that the curl of 
$P_\epsilon $ vanishes at the origin implies that $d\hat g_\epsilon 
(0) = 0$ for all $\epsilon $, see section 2. Then the 2-jet of 
$\hat g_\epsilon $ is $R$-equivalent to $(x^2\pm y^2)$, and there is a 
formal 
change of coordinates of the form (4.12) preserving the  
volume form $dx\wedge dy$ and reducing $\hat g_\epsilon $ to 
the form $$\lambda _1(\epsilon )(x^2\pm y^2) + 
\lambda _2(\epsilon )(x^2\pm y^2)^2 + 
\cdots . \eqno (4.13)$$ By the condition (1.4) 
this change of coordinates brings 
the series $\hat f_\epsilon $ to the form 
$$\mu _1(\epsilon )(x^2\pm y^2) + \mu _2(\epsilon )(x^2\pm y^2)^2 + 
\cdots , \eqno (4.14)$$ 
where $\mu _1(0) = 0$ (the latter 
follows from the definiton of $N$ singularities).    
So, $P_\epsilon $ reduces to the formal normal form (1.3), 
where the formal  
series $\hat g_\epsilon $ and $\hat f_\epsilon $ have the form 
(4.13) and (4.14). Now we can return to the smooth normal form 
(4.8), where the formal series of $c_\epsilon, d_\epsilon , 
a_\epsilon $ and $b_\epsilon $ are equal, respectively, to  
${\partial \hat f_\epsilon \over \partial x},  
{\partial \hat f_\epsilon \over \partial y}, 
{\partial \hat g_\epsilon \over \partial x}, 
{\partial \hat g_\epsilon \over \partial y}$. Proposition 4.1  
is a direct corollary of this normal form and the relation (4.10). 

\v 

\centerline {\bf 5. Proof of Theorem 1.1.} 

\v

It is easy to see that the condition that 
$j_0^1P$ is not isomorphic to (1.1) implies that there are three 
vanishing at the origin differentially independent 
functions $a,b$ and $c$  such that 
$P(da,db) = c$. Then $P_\epsilon (da,db) = c_\epsilon $, where 
$c_\epsilon $ is a family of functions such that $c_0 = c$, and in the 
coordinate system $x=a, y=b, z = c_\epsilon $ the family $P_\epsilon $ 
has the form (1.2). 

In what follows we will use the following notation. By 
$a^{(i)}_\epsilon (x,y,z), b^{(i)}_\epsilon (x,y,z), \dots $ 
we denote functions which are 
homogeneous degree $i$ polynomials with respect to the first two 
coordinates with coefficients being smooth functions of $\epsilon $ 
and the last coordinate. Also, we will say that a function 
$h_\epsilon (x,y,z)$ is affine with respect to $z$ if it can be 
written in the form $h_{0,\epsilon }(x,y) + zh_{1,\epsilon }(x,y)$. 
Let us show that the normal form (1.3)-(1.4) 
is a corollary of the following lemma. 

\v 

{\bf Lemma 5.1.} {\it Let $q\ge 0$ and $P_\epsilon $ be a family of 
Poisson structures of the form 

$$\{x,y\} = z + c^{(q+2)}_\epsilon (x,y,z) + c^{(q+3)}_\epsilon (x,y,z) + 
\cdots, $$ 

$$ \{y,z\} = a^{(0)}_\epsilon (x,y,z) + a^{(1)}_\epsilon (x,y,z) + 
\cdots + a^{(q-1)}_\epsilon (x,y,z) + a^{(q)}_\epsilon (x,y,z) + \cdots , 
\eqno (5.1)$$ 

$$ \{z,x\} = b^{(0)}_\epsilon (x,y,z) + b^{(1)}_\epsilon (x,y,z) + 
\cdots + b^{(q-1)}_\epsilon (x,y,z) + b^{(q)}_\epsilon (x,y,z) + \cdots , $$ 

\noindent where the functions $a^{(0)}_\epsilon , \dots , a^{(q-1)}_\epsilon $ and 
$b^{(0)}_\epsilon , \dots , b^{(q-1)}_\epsilon $ are affine with respect 
to $z$. There exists a change of coordinates of the form 
$$X = x + \mu _\epsilon ^{(q+2)}(x,y,z), \ Y = y, \ 
Z = z + \gamma _\epsilon ^{(q+1)}(x,y,z) +  
\gamma _\epsilon ^{(q+2)}(x,y,z)$$ 
reducing $P_\epsilon $ to the form 

$$\{X,Y\} = Z + C^{(q+3)}_\epsilon 
(X,Y,Z) + 
C^{(q+4)}_\epsilon (X,Y,Z) + 
\cdots, $$ 

$$ \{Y,Z\} = a^{(0)}_\epsilon (X,Y,Z)  
+ a^{(1)}_\epsilon (X,Y,Z) + 
\cdots + a^{(q-1)}_\epsilon (X,Y,Z) + 
A^{(q)}_\epsilon (X,Y,Z) + \cdots , 
\eqno (5.2)$$ 

$$ \{Z,X\} = b^{(0)}_\epsilon (X,Y,Z)  
+ b^{(1)}_\epsilon (X,Y,Z) +  
\cdots + b^{(q-1)}_\epsilon (X,Y,Z) + 
B^{(q)}_\epsilon (X,Y,Z) + 
\cdots , $$ 

\noindent where the functions $A_\epsilon ^{(q)}$ and 
$B_\epsilon ^{(q)}$ are affine with respect to $Z$.}

\v

This lemma implies the existence of a family of formal change of 
coordinates (centered at the origin for all $\epsilon $) reducing 
any Poisson structure of the form (1.2) to the form 
$$\{x,y\} = z, \ \{y,z\} = a_\epsilon (x,y) + zb_\epsilon (x,y), \ 
\{y,z\} = c_\epsilon (x,y) + zd_\epsilon (x,y),$$ with some formal 
series $a_\epsilon , b_\epsilon , c_\epsilon $ and $d_\epsilon $. 
The Jacobi identity gives the relations 
$$a_\epsilon d_\epsilon = b_\epsilon d_\epsilon , \ 
{\partial a_\epsilon \over \partial x} = 
{\partial c_\epsilon \over \partial y}, \   
{\partial b_\epsilon \over \partial x} = 
{\partial d_\epsilon \over \partial y},$$ 
and the normal form (1.3)-(1.4) follows. 

\v 

{\bf Proof of Lemma 5.1.} Make a change of coordinates of the form 
$$\tilde x = x, \ \tilde y = y, \ \tilde z = z(1 + e_\epsilon ^{
(q+1)}(x,y,z)).\eqno (5.3)$$ This change transforms (5.1) to the form 

$$\{\tilde x,\tilde y\} = \tilde z + z\tilde c^{(q+1)}_\epsilon 
(\tilde x,\tilde y,\tilde z) + 
\tilde c^{(q+2)}_\epsilon (\tilde x,\tilde y,\tilde z) + 
\cdots, $$ 

$$ \{\tilde y,\tilde z\} = a^{(0)}_\epsilon (\tilde x,\tilde y,\tilde z) 
+ a^{(1)}_\epsilon (\tilde x,\tilde y,\tilde z) + 
\cdots + a^{(q-1)}_\epsilon (\tilde x,\tilde y,\tilde z) + 
A^{(q)}_\epsilon (\tilde x,\tilde y,\tilde z) + \cdots , 
\eqno (5.4)$$ 

$$ \{\tilde z,\tilde x\} = b^{(0)}_\epsilon (\tilde x,\tilde y,\tilde z) 
+ b^{(1)}_\epsilon (\tilde x,\tilde y,\tilde z) + 
\cdots + b^{(q-1)}_\epsilon (\tilde x,\tilde y,\tilde z) + 
B^{(q)}_\epsilon (\tilde x,\tilde y,\tilde z) + 
\cdots , $$   

\noindent where 
$$A_\epsilon ^{(q)} = a_\epsilon ^{(q)} - 
z^2{\partial e_\epsilon ^{(q+1)}\over \partial x}, \  
B_\epsilon ^{(q)} = b_\epsilon ^{(q)} - 
z^2{\partial e_\epsilon ^{(q+1)}\over \partial y}.\eqno (5.5)$$ 
Let 

$$A_\epsilon ^{(q)} = \alpha _{0, \epsilon }(x,y) +  
z\alpha _{1, \epsilon }(x,y) + z^2 \alpha _{2, \epsilon }(x,y,z), $$  

$$B_\epsilon ^{(q)} = \beta _{0, \epsilon }(x,y) +  
z\beta _{1, \epsilon }(x,y) + z^2\beta _{2, \epsilon }(x,y,z).$$ 

\noindent Write the Jacobi identity for the Poisson structure (5.1)   
in the form  
$$E_{0,\epsilon }(x,y) + zE_{1,\epsilon }(x,y) + 
z^2E_{2,\epsilon }(x,y) + z^3E_{3,\epsilon }(x,y,z) = 0.$$ 
Using the condition that 
the functions $a^{(0)}_\epsilon , \dots , a^{(q-1)}_\epsilon $ and 
$b^{(0)}_\epsilon , \dots , b^{(q-1)}_\epsilon $ are affine with respect 
to $z$, we obtain 
$$E_{3,\epsilon }(x,y,z) = E_{3,\epsilon }^{(q-1)}(x,y,z) + 
E_{3,\epsilon }^{(q-2)}(x,y,z) + \cdots , \ \ \ 
E_{3,\epsilon }^{(q-1)}(x,y,z) = z^3({\partial \alpha _{2,\epsilon }
\over \partial y} - {\partial \beta _{2,\epsilon }
\over \partial x}).$$ Therefore  
$${\partial \alpha _{2,\epsilon }
\over \partial y} - {\partial \beta _{2,\epsilon }
\over \partial x} = 0.\eqno (5.6)$$  
The relations (5.5) and (5.6) imply that there exists a change 
of coordinates of the form (5.3) which reduces (5.1) to the form 
(5.4), where the functions $A_\epsilon ^{(q)}$ and 
$B_\epsilon ^{(q)}$ are affine with respect to $\tilde z$. 

Our next step is a change of coordinates of the form 
$$\hat x = \tilde x + r_\epsilon ^{
(q+2)}(\tilde x,\tilde y,\tilde z), \ \hat y = \tilde y, \ \hat z = 
\tilde z,\eqno (5.7)$$ 
where 
$${\partial r_\epsilon ^{(q+2)}\over \partial x} = - \tilde c_\epsilon 
^{(q+1)}.$$ It is easy to see that the change (5.7) reduces (5.4) to  
the form 

$$\{\hat x,\hat y\} = \hat z  + 
\hat c^{(q+2)}_\epsilon (\hat x,\hat y,\hat z) + 
\hat c^{(q+3)}_\epsilon (\hat x,\hat y,\hat z) +  
\cdots, $$ 

$$ \{\hat y,\hat z\} = a^{(0)}_\epsilon (\hat x,\hat y,\hat z) 
+ a^{(1)}_\epsilon (\hat x,\hat y,\hat z) + 
\cdots + a^{(q-1)}_\epsilon (\hat x,\hat y,\hat z) + 
A^{(q)}_\epsilon (\hat x,\hat y,\hat z) + \cdots , 
\eqno (5.8)$$ 

$$ \{\hat z,\hat x\} = b^{(0)}_\epsilon (\hat x,\hat y,\hat z) 
+ b^{(1)}_\epsilon (\hat x,\hat y,\hat z) + 
\cdots + b^{(q-1)}_\epsilon (\hat x,\hat y,\hat z) + 
B^{(q)}_\epsilon (\hat x,\hat y,\hat z) + 
\cdots . $$ 

Finally, to reduce (5.8) to the requirted normal form 
(5.2) it suffices to make a change of coordinates 
$$X = \hat x, \  Y = \hat y, \ Z = \hat z  + 
\hat c^{(q+2)}_\epsilon (\hat x,\hat y,\hat z).$$

\v 

\centerline {\bf 6. Appendix. Normal form for integrable 1-forms.}

\v 

Let $P$ be a local Poisson structure on $R^3$, and let $\Omega $ be a 
local nondegenerate volume form. Consider $P$ as a field of 2-vectors, 
then we can associate to $P$ a differential 1-form 
$\omega $ such that $\omega (Y) = \Omega (Y\wedge P)$ for any vector field 
$Y$. The relation $[P,P] = 0$, valid for any Poisson structure, implies 
the integrability of $\omega $: $\omega \wedge d\omega = 0$. The 
integrable 1-form $\omega $ depends on the choice of $\Omega $, and it is 
clear that $\omega $ is invariantly related to $P$ up to multiplication 
by a nonvanishing function. Denote by $(\omega )$ the Pfaffian equation 
generated by $\omega $, i.e., a module of differential 1-forms over the 
ring of smooth functions generated by $\omega $. This 
Pfaffian equation is invariantly related to the Poisson structure $P$, 
and we will denote it by $b(P)$. If, in local coordinates, $P$ has the 
form 

$$\{x,y\} = B(x,y,z), \  \{y,z\} = C(x,y,z), \  \{z,x\} = D(x,y,z),$$ 

\noindent then the Pfaffian equation $b(P)$ is generated by the 1-form 

$$B(x,y,z)dz + C(x,y,z)dx + D(x,y,z)dy.$$ 

Note that any local integrable Pfaffian equation $(\omega )$ 
has the form $b(P)$ 
for some Poisson structure $P$ (invariantly related to $(\omega )$ up 
to multiplication by a nonvanishing function), 
therefore all 
the results of this paper imply immediate corollaries 
on the classification 
of integrable Pfaffian equations. In particular, the corollary of 
our main result, Theorem 1.1, is as follows. 

\v

{\bf Theorem 6.1.} {\it Let $\omega _\epsilon $ be a local family of 
integrable 1-forms on $R^3$ such that $\omega _0(0) = 0$ and 
$j^1_0\omega _0$ is not equivalent to $a(xdy - ydx)$, $a\in R$. 
Then the family of Pfaffian equations $(\omega _\epsilon )$ is  
equivalent to a family of the form $$(zdz + U_\epsilon (x,y,z)dx +  
V_\epsilon (x,y,z)dy),\eqno (6.1)$$ 
and there exists a family of formal changes of 
the coordinates $x,y,z$, 
centered at the point $x=y=z=0$ for all $\epsilon $,     
reducing (6.1) to a family of Pfaffian equations generated by 
1-forms $$zdz + d\hat f_\epsilon (x,y) + zd\hat g_\epsilon (x,y),
\eqno (6.2)$$ where the formal series $\hat f_\epsilon $ and 
$\hat g_\epsilon $ satisfy the relation $d\hat f_\epsilon \wedge 
d\hat g_\epsilon \equiv 0.$} 

\v 

The singularities of integrable 1-forms were studies by many authors.  
The $V$ singularities of Poisson structures 
correspond to the case where $d\hat g_0(0)\ne 0$, 
i.e., to integrable Pfaffian equation  generated by a 1-form 
$\omega $ such that $d\omega (0)\ne 0$. In this case one can use the 
Darboux theorem on classification of closed nondegenerated 2-forms 
to show that $\omega $ is equivalent to a 1-form of the form 
$a(x,y)dx + b(x,y)dy$. Therefore the classification of $V$ 
singularities of integrable Pfaffian equations reduces to the orbital 
classification of vector fields on the plane. This well known 
reduction (see [Ku]) is similar to the results of section 2.1. 

\v 

The $sl(2), so(3), A^+$ and $A^-$ singularities of Poisson structures 
correspond, respectively, to the cases where in the normal form (6.2) 
the formal series $\hat f_0$ is $R$-equivalent to $x^2 + y^2, \ 
x^2 - y^2, \ x^2\pm y^m, \ -x^2\pm y^m$, \ $m\ge 3$, and $\hat g_0$ has 
the form $\hat g_0 = \hat \lambda \circ \hat f_0$ for some formal 
series $\hat \lambda $. These singularities 
of integrable Pfaffian equations are algebraically isolated. The results 
of the papers [Mo] and [Ma] imply that they can be reduced to a formal 
normal form $(dW)$, where $W = z^2 + x^2\pm y^2$ in the $so(3)$ and 
$sl(2)$ cases, and $W = z^2\pm x^2\pm y^m$ in the $A$ case. These 
normal forms also follow from the works [Co1], [Co2], [We1] and 
our results in section 3.   

\v 

Using the 1-1 correspondence between Poisson structures and integrable 
Pfaffian equations, we also can make corollaries of our results 
concerning the  
$N$ singularities of integrable Pfaffian equations (the case where 
$\hat g_0$ is $R$-equivalent to $x^2\pm y^2$ and $j_0^2\hat f_0 = 0$) 
and bifurcations of integrable Pfaffian 
equations near $A$ and $N$ singularities. As far as we know, 
these results and Theorem 6.1 are new.      

\v 

Theorem 6.1 can be generalized to the $n$-dimensional case. 
Such a generalization is an analogous of a generalization of 
Theorem 1.1 to Nambu structures. The 
corresponding results and their corollaries will be published 
elsewhere.

\v

\parindent=0mm

\centerline {\bf References} 

\v 

[AI] V.I. Arnold, Yu.S. Ilyashenko, Ordinary differential 
equations, {\it in Encyclopedia of Mathematical Sciences, 
Dynamical Systems - 1, Springer-Verlag, 1988.} 

[Ar] V.I. Arnold, Mathematical methods of classical mechanics, 2nd  
edition, {\it Reading, Benjamin and Cummings, 1978.}

[AVG] V.I. Arnold, A.N. Varchenko, S.M. Gusein-Zade, 
Singularities of differentiable maps - 1, {\it Birkhauzer, 1985.}   

[Bo] R.I. Bogdanov, Modules of $C^\infty $ orbital normal forms 
of singular points of vector fields on a plane, {\it Functional 
Anal. Appl. 11, 1 (1977), pp. 57-58.}

[Co1] J.F. Conn, Normal forms for analytic Poisson structures, {\it Ann. 
of Math., 119 (2), 1984, pp. 576 - 601.}

[Co2] J.F. Conn, Normal forms for smooth Poisson structures, {\it Ann. 
of Math., 121 (2), 1985, pp. 565 - 593.}  

[DH] J.P. Dufour, A. Haraki, Rotationnels et structures de 
Poisson quadratiques, {\it C.R. Acad. Sci. Paris 312 I (1991), 
pp. 137-140.} 

[Du] J.P. Dufour, Lin\'earisation de certaines structures de Poisson, 
{\it J. Differential Geometry 32 (1990), pp. 415 - 428.}

[Ku] I. Kupka, The singularities of integrable structurally 
stable Pfaffian forms, {\it Proc. Acad. Sci. USA 52 (1964), p. 1431.}

[Li] A. Lichnerowicz, Les vari\'et\'es de Poisson et leurs alg\`ebres 
de Lie associ\'ees, {\it J. Differential Geometry 12 (1977), pp. 253-300.} 

[Ma] B. Malgrange, Frobenius avec singularit\'es. 1. Codimension un, 
{\it Publ. IHES 46 (1976), pp. 163 - 173.}

[Mo] R. Moussu, Sur l'existence d'int\'egrales premi\`eres pour un 
germe de forme de Pfaff, {\it Ann. Inst. Fourier 26, 2 (1976), 
pp. 171 - 220.}

[We1] A. Weinstein, The local structure of Poisson manifolds, 
{\it J. Differential Geometry 18 (1983), pp. 523 - 557.} 

[We2] A. Weinstein, The modular automorphism group of a 
Poisson manifold, {\it University of California at Berkeley, preprint.}

\bye